\documentclass[11pt, twoside]{amsart}

\usepackage[utf8]{inputenc}
\usepackage[T1]{fontenc}
\usepackage{amssymb,amsmath,amstext, graphicx}

\newtheorem{theor}{Theorem}[section]
\newtheorem{lem}[theor]{Lemma}
\newtheorem{defin}[theor]{Definition}

\newtheorem{prop}[theor]{Proposition} 

\newtheorem{notation}[theor]{Notation}

\newtheorem{examps}[theor]{Examples}
\newtheorem{cor}[theor]{Corollary}
\newtheorem{rem}[theor]{Remark}

\newtheorem{fact}[theor]{Fact}

\numberwithin{equation}{section}

\newcommand{\tp}{\mathrm{tp}}

\newcommand{\acl}{\mathrm{acl}}
\newcommand{\dcl}{\mathrm{dcl}}

\newcommand{\es}{\emptyset}

\newcommand{\su}{\mathrm{SU}}

\newcommand{\nts}{\negthickspace}

\newcommand{\approxc}{\nts \approx}
\newcommand{\crd}{\mathrm{crd}}

\newcommand{\meq}{^{\mathrm{eq}}}

\newcommand{\mcA}{\mathcal{A}}
\newcommand{\mcB}{\mathcal{B}}
\newcommand{\mcC}{\mathcal{C}}

\newcommand{\mcG}{\mathcal{G}}

\newcommand{\mcK}{\mathcal{K}}
\newcommand{\mcL}{\mathcal{L}}
\newcommand{\mcM}{\mathcal{M}}
\newcommand{\mcN}{\mathcal{N}}

\newcommand{\mcR}{\mathcal{R}}

\newcommand{\nind}{\raisebox{-2pt}[5pt][0pt]{$\smile$} 
\hspace*{-6.8pt}\raisebox{3pt}[5pt][0pt]{$|$}\hspace*{-6.8pt}
\raisebox{3pt}[5pt][0pt]{$\diagup$} }

\newcommand{\rng}{\mathrm{rng}}

\title[Homogeneous 1-based structures]
{Homogeneous 1-based structures and interpretability in random structures}

\author{Vera Koponen}

\address{Vera Koponen, Department of Mathematics, Uppsala University, Box 480,
75106 Uppsala, Sweden.}

\email{vera@math.uu.se}

\begin{document}

\begin{abstract}
Let $V$ be a finite relational vocabulary in which no symbol has arity greater than 2.
Let $\mcM$ be countable $V$-structure which is homogeneous, simple and 1-based.
The first main result says that if $\mcM$ is, in addition, primitive, then it is strongly interpretable in a random structure.
The second main result, which generalizes the first, implies (without the assumption on primitivity) that if $\mcM$ is 
``coordinatized'' by a set with SU-rank 1 and there is no definable (without parameters) nontrivial
equivalence relation on $M$ with only finite classes, then $\mcM$ is strongly interpretable in
a random structure.
\end{abstract}

\maketitle

\section{Introduction}

\noindent
A first-order structure $\mcM$ will be called {\em homogeneous} if it has a finite relational vocabulary 
and every isomorphism between finite substructures of $\mcM$ can be extended to an automorphism of $\mcM$.
For surveys of homogeneous structures and their connections to other areas see 
\cite{BP, Che98, HN, Mac10, Nes}.
Although there are $2^\omega$ countable nonisomorphic homogeneous structures for a vocabulary with only one binary relation symbol
\cite{Hen72},
it has been shown that in several cases, such as partial orders, undirected graphs, directed graphs and
stable relational structures (with finite relational  vocabulary), the countable homogeneous structures among them
can be classified in a rather concrete way \cite{Che98, Gar, GK, JTS, Lach84, Lach97, LT, LW, Schm, Shee}.
The work on stable homogeneous structures together and on ``geometric stability theory'' (notably by Zilber \cite{Zil})
was useful to reach a good understanding of
$\omega$-categorical $\omega$-stable structures \cite{CHL} and later of smoothly approximable structures \cite{CH, KLM}.

Stability theory was from the mid 90ies generalized to simplicity theory, where the class of simple theories is the largest class of
complete theories $T$ such that in every $\mcM \models T$ there is a symmetric ``independence relation'' on subsets of $M$.
The random graph (or Rado graph) is a standard example of a homogeneous simple structure (i.e. one which has simple theory)
which is not stable, and the same is true of
the ``random structure'' with respect to any finite relational vocabulary.
Since the {\em infinite} countable stable homogeneous structures are classified \cite{Lach97}, 
one may ask if it is possible to reach, if not a classification, at least
some systematic understanding of (infinite countable) simple homogeneous structures.
Besides the present work, \cite{AK14, Kop15} and the dissertation of Aranda L\'{o}pes \cite{AL} 
the author is not aware of any results in this direction.
The class of all simple homogeneous structures seems too wide to start with, 
so we focus on a subclass of it which seems easier to deal with.
We consider homogeneous, simple and 1-based structures, where the property {\em 1-based} implies 
(for structures in general) that the 
independence relation behaves ``nicely''
(like in a vector space or random structure, for example). 
In fact, it follows from work of Macpherson \cite{Mac91} and De Piro and Kim \cite{PK} that 
every homogeneous, simple and 1-based structure has trivial dependence.
If a homogeneous structure has only unary and binary relation symbols,
then it is 1-based if and only if it has trivial dependence
(see Definition~\ref{definition of trivial dependence} and 
Fact~\ref{homogeneity and 1-basedness implies finite rank and trivial dependence} below).

Moreover, we focus on binary homogeneous structures, where a structure is called {\em binary} if its vocabulary is finite, relational
and has no symbol with arity greater than~2. 
This simplifies the situation because the ``independence theorem'' of simple theories has strong consequences in the context
of binary structures. For example it follows, rougly speaking, 
that every $\es$-definable set of SU-rank 1 (in the extension by imaginaries)
is -- when viewed as a ``canonically embedded structure'' -- a reduct of a binary random structure.
The working hypothesis is that countable, binary, homogeneous, simple and 1-based structures are sufficiently uncomplicated that it should be possible to work out some sort of rather explicit understanding of them.
Moreover, the author does not know of any example of a homogeneous and simple structure which is not 1-based.
Such an example, particularly a binary one, would probably shed light on the understanding of simple homogeneous structures.
A proof that it does not exist would also be striking.

We say that a structure is {\em primitive} if there is no nontrivial equivalence relation on its universe which is definable without parameters.
(Nontrivial means that there are at least two classes and at least one has more than one element.)
It is easy to see that if $\mcM$ is homogeneous (and simple and 1-based) and nonprimitive, 
then the substructure on any equivalence class of the equivalence
relation that witnesses nonprimitivity is homogeneous (and simple and 1-based). 
So a good understanding of homogeneous simple structures requires an equally good understanding of
primitive homogeneous simple structures. In particular we have the following problem, where the notion of
a `binary random structure' is explained in Section~\ref{Random structures}:

\begin{itemize}
\item[] {\em Is every countable, binary, homogeneous, primitive, simple and 1-based structure a random structure?}
\end{itemize}

\noindent
The problem is open also if we remove the assumptions about binarity and 1-basedness.
Our first main result can be seen as an ``approximate solution'':

\begin{theor}\label{first main result}
Suppose that $\mcM$ is a structure which is countable, binary, homogeneous, primitive, simple and 1-based.
Then $\mcM$ is strongly interpretable in a  binary random structure.
\end{theor}

\noindent
The notion {\em strongly interpretable} (Definition~\ref{definition of strong interpretability}) 
implies the notion {\em interpretable} in its usual senses (e.g. \cite{Hod, Mac91}).
Hence Theorem~\ref{first main result} can be interpreted as saying that 
a countable, binary, homogeneous, primitive, simple and 1-based cannot be more complicated than 
a binary random structure. The theorem may also be of help to solve the problem stated above.

The second main result is a generalisation of the first one, but has a more complicated formulation, which
is why we stated Theorem~\ref{first main result} first and not as a corollary of Theorem~\ref{second main result} below.
Theorem~\ref{second main result} applies to some (but not all) {\em non}primitive binary homogeneous simple and 1-based structures.
For example, let $\mcR$ be the random graph, let $M = R^2$ and
let $\mcM$ be the structure with universe $M$ that has a binary relation symbol $R_p$ for every 4-type over $\es$ of $\mcR$
and where $R_p$ is interpreted as $\{((a,b), (a',b')) : \mcR \models p(a,b,a',b')\}$. 
Then $\mcM$ satisfies the hypotheses of Theorem~\ref{second main result}.
(The same holds if $\mcR$ is any binary random structure and we consider ordered $k$-tuples for any fixed $k > 1$.)
On the other hand if, for example, $\mcM$ is the structure with infinite countable universe,
two equivalence relations $R_1, R_2$, such that $R_2$ refines $R_1$ by splitting each $R_1$-class into infinitely many
infinite parts, then $\mcM$ does not have height 1, so Theorem~\ref{second main result} does not apply.\footnote{
This structure has height 2 in the sense of \cite{Djo06}. It is of course a structure which is easy to understand. 
But examples with height~2 and the other properties of $\mcM$ can be more complicated.}
The concept of `height', appearing below, is given by Definition~\ref{definition of height},
and the notation `$\mcM / \approx$'
is explained in Definition~\ref{definition of canonically embedded structure}:

\begin{theor}\label{second main result}
Suppose that $\mcM$ is a structure which is countable, binary, homogeneous, simple and 1-based with height 1.
Let $\approx$ denote the equivalence relation `$\acl_\mcM(x) = \acl_\mcM(y)$', where `$\acl_\mcM$' denotes algebraic closure in $\mcM$.\\
(i) Then there is a binary random structure $\mcR$ such that $\mcM / \approxc$ is definable in $\mcR\meq$.\\
(ii) If, in addition, the equivalence relation `$\approx$' is trivial, 
then $\mcM$ is strongly interpretable in a binary random structure.
\end{theor}

\noindent
The author hopes that this theorem can be generalised in a suitable way to `height $r$' for any $0 < r < \omega$.
One would then have a theorem which applies to all countable, binary, homogeneous, simple and 1-based structures
(since `1-based' implies `finite height'). 
But besides finding the ``right'' generalisation it seems like one has to overcome a number of technical difficulties.

This article is organised as follows.
Section~\ref{Preliminaries} recalls the necessary background about homogeneous structures and $\omega$-categorical simple structures.
The proofs of the above theorems are given in Section~\ref{Binary homogeneous 1-based structures},
which ends with a recipe for making structures that satisfy the main results of this article.
The proofs of the main results use the ``coordinatization results'' from \cite{Djo06},
which do not need the assumption that the structure is binary.
One also needs to know what 
a canonically embedded structure $\mcG$ in $\mcM\meq$ may look like if its universe is
a definable set with $\su$-rank 1. It turns out that if $\mcM$ is countable, binary, homogeneous and simple with trivial
dependence then $\mcG$ is a {\em reduct} (Definition~\ref{definition of reduct}) of a binary random structure. 
This is proved in \cite{AK14} and the result is refered to as Fact~\ref{rank 1 sets are reducts of binary random structures} below.
The main technical result (besides Fact~\ref{rank 1 sets are reducts of binary random structures}) is
Lemma~\ref{crd(a) is rigid} which is proved in Section~\ref{Proof of the main technical result}.
The assumption that $\mcM$ in Theorems~\ref{first main result} and~\ref{second main result}
is binary is only used in the application of Fact~\ref{rank 1 sets are reducts of binary random structures}
and in the proof of Lemma~\ref{crd(a) is rigid}.

\section{Preliminaries}\label{Preliminaries}

\noindent
The prerequisites of this article are more or less the same as those of \cite{AK14}.
We list, in this section, a number of definitions and facts, concerning homogeneous structures, $\omega$-categorical simple
structures and imaginary elements, in order to make this article relatively self contained,
but refer to \cite{AK14} for further explanations.

\subsection{General notation and terminology}

\noindent
We call a vocabulary (also called signature) {\em relational} if it only contains relation symbols.
Given a finite relational vocabulary the {\em maximal arity} of it is the largest integer $k$ such that some relation symbol of it has arity $k$.
If  $V$ is a {\em  finite} vocabulary and the maximal arity is 2 then we call $V$ {\em binary}
(although it may contain unary relation symbols), and in this case a $V$-structure may be called a {\em binary structure}.
We denote (first-order) structures by $\mcA, \mcB, \ldots, \mcM, \mcN, \ldots$ and their respective universes by 
$A, B, \ldots, M, N, \ldots$. 
Finite sequences (tuples) of elements of some structure (or set in general) will be denoted
$\bar{a}, \bar{b}, \ldots$, while $a, b, \ldots$ ususally denote elements from the universe of some structure.
The notation $\bar{a} \in A$ means that every element in the sequence $\bar{a}$ belongs to $A$. 
Sometimes we write $\bar{a} \in A^n$ to show that the length of $\bar{a}$, denoted $|\bar{a}|$, is $n$.
The {\em range of $\bar{a}$}, denoted $\rng(\bar{a})$, is the set of elements that occur in $\bar{a}$.
Notation regarding types, algebraic closure and definable closure is standard, where we may use
a subscript or superscript to indicate which structure we have in mind.
A structure $\mcM$ is called {\em $\omega$-categorical}, respectively {\em simple}, if 
its complete theory $Th(\mcM)$ has that property. (See \cite{Hod} and \cite{Cas, Wag} for definitions.)



\subsection{Homogeneous and $\omega$-categorical structures}
\label{Homogeneous structures}

\begin{defin}\label{definition of homogeneous}{\rm
(i) Let $V$ be a finite relational vocabulary and $\mcM$ a $V$-structure. We call $\mcM$ {\em homogeneous}
if for all finite substructures $\mcA$ and $\mcB$ of $\mcM$, every isomorphism from $\mcA$ to $\mcB$ can be extended to an
automorphism of $\mcM$.\\
(ii) We say that a structure $\mcM$, with any vocabulary, is {\em $\omega$-homogeneous} 
if whenever $n < \omega$,
$a_1, \ldots, a_n, a_{n+1}, b_1, \ldots, b_n \in M$ and $\tp_\mcM(a_1, \ldots, a_n) = \tp_\mcM(b_1, \ldots, b_n)$,
then there is $b_{n+1} \in M$ such that $\tp_\mcM(a_1, \ldots, a_{n+1}) = \tp_\mcM(b_1, \ldots, b_{n+1})$.
}\end{defin}

\noindent
For standard results about homogeneous structures we refer to (for example)~\cite[Sections~7.1 and~7.4]{Hod}.
We mainly use the following implications of homogeneity
(or $\omega$-categoricity), which follow from Corollary~7.4.2 in \cite{Hod} and from
the well known characterization of $\omega$-categorical theories by Engeler, Ryll-Nardzewski and Svenonius:

\begin{fact}\label{fact about homogeneous structures}
(i) If $\mcM$ is $\omega$-categorical then it is $\omega$-homogeneous.\\
(ii) Suppose that $\mcM$ is countable and $\omega$-homogeneous.
Then for all $0 < n < \omega$ and all $a_1, \ldots, a_n, b_1, \ldots, b_n \in M$ such that 
$\tp_\mcM(a_1, \ldots, a_n) = \tp_\mcM(b_1, \ldots, b_n)$, there is an automorphism $f$ of $\mcM$ such that
$f(a_i) = b_i$ for every $i$.\\
(iii) Suppose that $\mcM$ is a countable $V$-structure where $V$ is a finite relational vocabulary.
Then $\mcM$ is homogeneous if and only if $\mcM$ is $\omega$-categorical and has elimination of quantifiers.
\end{fact}

\subsection{Simple structures, imaginary elements, 1-basedness and triviality of dependence}\label{Simple omega-categorical structures}

\noindent
We assume familiarity with imaginary elements and $\mcM\meq$, defined in \cite{Hod, She} for example,
and with basic simplicity theory, as found in \cite{Cas, Wag} for example.
Since the distinction between sorts (of imaginary elements) will be relevant here, and since 
some notions and results are simplified when considering $\omega$-categorical simple theories, 
compared with simple theories in general, we will nevertheless rehearse some notions and results that will be used.

Let $V$ be a vocabulary and $\mcM$ a $V$-structure.
For every $0 < n < \omega$ and $\es$-definable equivalence relation $E$ on $M^n$,
$V\meq$ (the vocabulary of $\mcM\meq$) contains a unary relation symbol $P_E$ and 
a relation symbol $F_E$ of arity $n+1$ (both of which do not belong to $V$), where $P_E$ is
interpreted in $\mcM\meq$ as the set of $E$-equivalence classes and
$F_E$ is interpreted as the graph of the function which sends every $\bar{a} \in  M^n$
to its equivalence class.
A {\em sort} of $\mcM\meq$ is, by definition, a set of the form $S_E = \{a \in M\meq : \mcM\meq \models P_E(a) \}$ for some 
$E$ as above. If $A \subseteq M\meq$ and there are only finitely many $E$ such that $A \cap S_E \neq \es$, 
then we say that {\em only finitely many sorts are represented in $A$}.
The identity relation, `=', is clearly a $\es$-definable equivalence relation on $M$ and every $=$-class is a singleton.
Therefore $M$ can (and will) be identified with the sort $S_=$, which we call the {\em real sort}, so $M \subseteq M\meq$.
Below follow some facts and definitions. See for example \cite{AK14} for explanations or proofs of these facts.

\begin{fact}\label{facts about simple omega-categorical structures}
Suppose that $\mcM$ is $\omega$-categorical and countable, 
let  $A \subseteq M\meq$ and suppose that only finitely many sorts are represented in $A$. \\
(i) For every $n < \omega$ and finite $B \subseteq \mcM\meq$, only finitely many types from $S_n^{\mcM\meq}(\acl_{\mcM\meq}(B))$
are realized by $n$-tuples in $A^n$.\\
(ii) If $B \subseteq M\meq$ is finite and $\bar{a} \in M\meq$,
then $\tp_{\mcM\meq}(\bar{a} / \acl_{\mcM\meq}(B))$ is isolated. \\
(iii) If $B \subseteq \mcM\meq$ is finite, $n < \omega$ and $p \in S_n^{\mcM\meq}(\acl_{\mcM\meq}(B))$ is realized 
in some elementary extension of $\mcM\meq$ by an $n$-tuple of imaginary elements
(i.e. elements satisfying $P_E(x)$ for some, not necessarily the same, $E$),
then $p$ is realized in $\mcM\meq$.
\end{fact}

\begin{defin}\label{definition of canonically embedded structure}{\rm
(i) We say that a structure $\mcN$ is {\em canonically embedded} in $\mcM\meq$ if $N$ is a $\es$-definable subset of $M\meq$ and
for every $0 < n < \omega$ and every relation $R \subseteq N^n$ which is $\es$-definable in $\mcM\meq$ there is a
relation symbol in the vocabulary of $\mcN$ which is interpreted as $R$ and the vocabulary of $\mcN$ contains no other relation symbols
(and no constant or function symbols).\\
(ii) In particular, if $R$ is an equivalence relation on $M$ which is $\es$-definable in $\mcM$,
then $\mcM/R$ denotes the canonically embedded structure with universe $M/R$, 
where $M/R$ is the set of all equivalence classes of $R$ (which is a $\es$-definable subset of $M\meq$).
}\end{defin}

\noindent
We immediately get the following:

\begin{fact}\label{first fact about canonically embedded structures}
If $\mcN$ is canonically embedded in $\mcM\meq$, then for all $\bar{a}, \bar{b} \in N$ and all $C \subseteq N$,
$\acl_\mcN(C) = \acl_{\mcM\meq}(C) \cap N$ and
$\tp_\mcN(\bar{a} / C) = \tp_\mcN(\bar{b} / C)$ if and only if $\tp_{\mcM\meq}(\bar{a} / C) = \tp_{\mcM\meq}(\bar{b} / C)$. 
\end{fact}

\noindent
Suppose that $T$ is a complete simple theory. 
For every type $p$ (possibly over a set of parameters) with respect to $T$, there is a notion of SU-rank of $p$, denoted $\su(p)$
(which is either an ordinal or undefined);

See for instance~\cite{Cas, Wag} for definitions and basic results about SU-rank. 
As usual, we abbreviate $\su(\tp_\mcM(\bar{a} / A))$ with $\su(\bar{a} / A)$
and $\su(\bar{a} / \es)$ with $\su(\bar{a})$.
(When using this notation there will be no ambiguity about which structure we work in.)
If $\su(\bar{a})$ is finite for every $\mcM \models T$ and every $\bar{a} \in M$, then we say that $T$ (and any $\mcM \models T$) has {\em finite $\su$-rank}.

\begin{defin}\label{definition of trivial dependence}{\rm
Let $T$ be a complete simple theory.\\
(i) We say that  $T$ has {\em trivial dependence} if, whenever $\mcM \models T$,
$A, B, C_1, C_2 \subseteq M\meq$ and $A \underset{B}{\nind} (C_1 \cup C_2)$, 
then $A \underset{B}{\nind} C_i$ for $i = 1$ or $i = 2$.\\
(ii) $T$ (as well as every model of it) is {\em 1-based} if for every $\mcM \models T$ and all $A, B \subseteq M\meq$, 
$A$ is independent from $B$ over $\acl_{\mcM\meq}(A) \cap \acl_{\mcM\meq}(B)$.
}\end{defin}

\noindent
For homogeneous structures the notions of 1-basedness, triviality of dependence and finiteness of rank 
are fairly tightly connected, in particular in the binary case.

\begin{fact}\label{homogeneity and 1-basedness implies finite rank and trivial dependence}
(i) Suppose that $\mcM$ is homogeneous, simple and 1-based.
Then $Th(\mcM)$ has trivial dependence and finite $\su$-rank (so in particular it is supersimple).\\
(ii) Suppose that $\mcM$ is binary, homogeneous and simple.
Then $\mcM$ has finite SU-rank. Also, the following three conditions are equivalent:
(a) $\mcM$ is 1-based, (b) $\mcM$ has trivial dependence, and 
(c) every type (over a finite set of parameters) of SU-rank 1 has trivial pregeometry (given by algebraic closure restricted
to that type).
\end{fact}

\noindent
{\bf Proof.}
(i) Let $\mcM$ satisfy the premisses of part~(i) of the lemma, so it follows that $\mcM$ is $\omega$-categorical.
By Theorem~1.1 in \cite{Mac91}, it is not possible to interpret an infinite group in $\mcM$, which, with the
terminology of \cite{Mac91}, means that it is not possible to define, with finitely many parameters, an infinite group in $\mcM\meq$.
Corollary~3.23 in \cite{PK} implies that if $Th(\mcM)$ does not have trivial dependence,
then an infinite group is definable, with finitely many parameters, in $\mcM\meq$.
It follows that $\mcM$ must have trivial dependence.
And finally, Corollary~4.7 in \cite{HKP} says that every simple, 1-based and $\omega$-categorical theory is supersimple with 
finite $\su$-rank.

(ii) The first claim is the main result of~\cite{Kop15}.
The second claim follows by combining results from~\cite{PK, HKP, Mac91} and is explained in some more detail
in the introduction to~\cite{Kop15}.
\hfill $\square$

\begin{defin}\label{definition of height}{\rm
Let $\mcM$ be a simple structure.
We say that $\mcM$ has {\em height 1} if there is a $\es$-definable $D \subseteq M\meq$ in which only finitely 
many sorts are represented and  $M \subseteq \acl_{\mcM\meq}(D)$ and $\su(d) = 1$ for every $d \in D$.
}\end{defin}

\subsection{Random structures}\label{Random structures}

\begin{defin}\label{definition of a binary random structure}{\rm
Let $V$ be a binary vocabulary and let $\mcM$ be an infinite homogeneous $V$-structure.\\
(i) A finite $V$-structure $\mcA$ is called a {\em forbidden structure with respect to $\mcM$}
if $\mcA$ cannot be embedded into $\mcM$.\\
(ii) Suppose that $\mcA$ is a forbidden structure with respect to $\mcM$.
We call $\mcA$ a {\em minimal} forbidden structure with respect to $\mcM$ if 
no proper substructure of $\mcA$ is a forbidden structure with respect to $\mcM$.\\
(iii) We call {\em $\mcM$ a binary random structure} if there does {\em not} exist a minimal forbidden structure $\mcA$ with
respect to $\mcM$ such that $|A| \geq 3$.
}\end{defin}

\begin{rem}\label{remarks about binary random structures}{\rm
(i) It is straightforward to see that the above definition of a binary random structure is equivalent
to the definition given in~\cite[Section~2.3]{AK14}.\\
(ii) The Rado graph, usually called random graph in model theory, is of course an example of a binary random structure.\\
(iii) Suppose that $\mcM$ is a binary random structure. Let $x_1, \ldots, x_n$, where $n \geq 3$, be distinct variables
and, for all $1 \leq i < j \leq n$, let $p_{i, j}(x_i, x_j) \in S^\mcM_2(\es)$.
Moreover assume that for all $i < j$ and all $i' < j'$, if $k \in \{i, j\} \cap \{i', j'\}$, then
the restriction of $p_{i,j}$ to $x_k$ is identical to the restriction of $p_{i',j'}$ to $x_k$.
It now follows straightforwardly from the definition of a binary random structure that 
$\bigcup_{1 \leq i < j \leq n} p_{i, j}(x_i, x_j)$ is consistent with $Th(\mcM)$ and realized in $\mcM$.
Since $\mcM$ is homogeneous it also follows that, for any $1 \leq i < j \leq n$, 
if $\mcM \models p_{i,j}(a_i, a_j)$ then there are
$a_k \in M$ for all $k \in \{1, \ldots, n\} \setminus \{i, j\}$ such that 
$\mcM \models \bigwedge_{1 \leq k < l \leq n} p_{k, l}(a_k, a_l)$.
}\end{rem}

\begin{fact}\label{properties of binary random structures}
Let $\mcM$ be a binary random structure.
Then $Th(\mcM)$ is simple, has $\su$-rank 1 and is 1-based with trivial dependence.
\end{fact}

\noindent
{\bf Proof sketch.}
Let $\mcM$ be a binary random structure, so it is homogeneous.
That $Th(\mcM)$ is simple, has $\su$-rank 1 and trivial dependence 
is proved in essentially the same way as the (folkore) result that the random graph has these properties.
It now follows from Corollary~4.7 in \cite{HKP} (where the terminology `modular'
is used in stead of '1-based') that $Th(\mcM)$ is 1-based.
\hfill $\square$
\\

\noindent
The proofs of the main results, Theorems~\ref{first main result} and~\ref{second main result}, 
use the notion of reduct and Fact~\ref{rank 1 sets are reducts of binary random structures}, below,
from~\cite{AK14}.

\begin{defin}\label{definition of reduct}{\rm
Let $\mcM$ and $\mcN$ be structures which need not have the same vocabulary.
We say that $\mcM$ is a {\em reduct} of $\mcN$ if they have the same universe ($M = N$) and
for every $0 < n < \omega$, if $R \subseteq M^n$ is $\es$-definable in $\mcM$, then it is $\es$-definable in $\mcN$.
}\end{defin}

\begin{fact}\label{rank 1 sets are reducts of binary random structures} {\rm \cite{AK14}}
Let $\mcM$ be countable, binary, homogeneous and simple with trivial dependence.
Suppose that $G \subseteq M\meq$ is $\es$-definable, only finitely many sorts are represented in $G$,
and $\su(a) = 1$ and $\acl_{\mcM\meq}(a) \cap G = \{a\}$ for every $a \in G$.
Let $\mcG$ denote the canonically embedded structure in $\mcM\meq$ with universe $G$.
Then $\mcG$ is a reduct of a binary random structure.
\end{fact}

\subsection{Interpretability}\label{Interpretability}

\begin{defin}\label{definition of strong interpretability}{\rm
Let $\mcM$ and $\mcN$ be structures, possibly with different vocabularies.\\
(i) We say that $\mcN$ is {\em strongly interpretable in $\mcM$} if there are
\begin{itemize}
\item $0 < n < \omega$, 
\item a formula $\chi(x_1, \ldots, x_n)$ without parameters in the language of $\mcM$,
\item a bijective function $f : \chi(\mcM) \to N$, and
\item for every $0 < k < \omega$ and formula $\varphi(x_1, \ldots, x_k)$ without parameters in the language of $\mcN$,
a formula $\psi_\varphi(\bar{y}_1, \ldots, \bar{y}_k)$ without parameters in the language of $\mcM$,
such that, for all $\bar{a}_1, \ldots, \bar{a}_k \in \chi(\mcM)$, 
\[ \mcM \models \psi_\varphi(\bar{a}_1, \ldots, \bar{a}_k) \ \Longleftrightarrow \ \mcN \models \varphi(f(\bar{a}_1), \ldots, f(\bar{a}_k)).\]
\end{itemize}
(ii) We say that $\mcN$ is {\em definable in $\mcM$} if $\mcN$ is strongly interpretable in $\mcM$ and it is possible to
choose $n = 1$ in the definition of `strongly interpretable'.
}\end{defin}

\noindent
It is immediate that if $\mcN$ is strongly interpretable in $\mcM$, 
then $\mcN$ is interpretable in $\mcM$ in the sense of Chapter~5.3 in \cite{Hod}, and $\mcN$ is interpretable in $\mcM$ 
in the sense of \cite{Mac91}, and it is definable in $\mcM\meq$ in the sense of \cite{PK}.
The following will be convenient to use.

\begin{lem}\label{equivalent condition of being strongly interpretable}
Suppose that there are positive integers $l, n_1, \ldots, n_l$, formulas $\chi_i(x_1, \ldots, x_{n_i})$ for $i = 1, \ldots, l$,
a bijective function $f : \bigcup_{i=1}^l \chi_i(\mcM) \to N$ and, for every
$0 < k < \omega$, $1 \leq i \leq l$ and formula $\varphi(x_1, \ldots, x_k)$ without parameters in the language of $\mcN$,
a formula $\psi_{\varphi,i}(\bar{y}_1, \ldots, \bar{y}_k)$ without parameters in the language of $\mcM$,
such that, for all $\bar{a}_1, \ldots, \bar{a}_k \in \chi_i(\mcM)$, 
\[ \mcM \models \psi_{\varphi,i}(\bar{a}_1, \ldots, \bar{a}_k) \ \Longleftrightarrow \ \mcN \models \varphi(f(\bar{a}_1), \ldots, f(\bar{a}_k)).\]
Then $\mcN$ is strongly interpretable in $\mcM$.
\end{lem}

\noindent
{\bf Proof sketch.}
Let $n$ be the maximum of $n_1, \ldots, n_l$.
Now the idea is that for every $i = 1, \ldots, l$ and $\bar{a} = (a_1, \ldots, a_i) \in \chi_i(\mcM)$, 
$\bar{a}$ can be ``translated'' into an $n$-tuple $\bar{a}' = (a'_1, \ldots, a'_n)$ where $a'_j = a_j$ for $j = 1, \ldots, n_i$
and $a'_j = a'_{n_i}$, for $j = n_i, \ldots, n$. In this way the set 
$\{ \bar{a}' : \bar{a} \in \chi_i(\mcM) \}$
is $\es$-definable in $\mcM$, and it follows that the union of these sets, for $i = 1, \ldots, l$, is also $\es$-definable in $\mcM$.
The rest is straightforward, via obvious modifications of $f$ and $\psi_{\varphi,i}$ for each $\varphi$.
\hfill $\square$

\section{Binary homogeneous 1-based structures}
\label{Binary homogeneous 1-based structures}

\noindent
In this section we prove the main results, Theorems~\ref{first main result} and~\ref{second main result}.
Throughout this section we make the following assumption, which is shared by both theorems:

\begin{itemize}
\item[]
{\bf $\mcM$ is countable, binary, homogeneous, simple and 1-based.}
\end{itemize}

\noindent
Now Fact~\ref{homogeneity and 1-basedness implies finite rank and trivial dependence} implies that

\begin{itemize}
\item[] {\bf $Th(\mcM)$ has trivial dependence and is supersimple with finite rank.}
\end{itemize}

\noindent
Throughout this section and Section~\ref{Proof of the main technical result} we use the following notational convention:

\begin{notation}\label{notation for acl etc}{\rm
The notations 
$\acl( \ )$, $\dcl( \ )$ and $\tp( \ )$ are abbreviations of
$\acl_{\mcM\meq}( \ )$, $\dcl_{\mcM\meq}( \ )$ and $\tp_{\mcM\meq}( \ )$, respectively.
However, when speaking of algebraic closure, definable closure or types in some other structure ($\mcM$ for example),
then we will show this explicitly with a subscript.
}\end{notation}

\noindent
The first step is is to use results from \cite{Djo06} to show that every $a \in M$ has finitely many ``coordinates'' of rank 1 in $\mcM\meq$
which to a large extent determine the properties of $a$.
By the results in Section~3 of \cite{Djo06},
there is a {\em self-coordinatized} set $C \subseteq M\meq$, in the sense of Definition 3.3 in \cite{Djo06}, such that
\begin{align}\label{C_i}
&\text{  $M \subseteq C$, $C$ is $\es$-definable;} \\
&\text{ there are $0 < r < \omega$ and $\es$-definable sets 
$C_0 \subseteq C_1 \subseteq \ldots \subseteq C_r \subseteq C$ such that} \nonumber \\
&\text{ $C_0 = \es$, and for every $n < r$ and every $a \in C_{n+1}$, $\su(a/C_n) = 1$;} \nonumber \\
&\text{ for all $n \leq r$, if $a \in C_n$, $b \in M\meq$ and $\tp(a) = \tp(b)$, then $b \in C_n$;} \nonumber \\
&\text{ only finitely many sorts are represented in $C$; and} \nonumber \\
&\text{ $M \subseteq C \subseteq \acl(C_r)$.} \nonumber
\end{align}

\noindent
We assume that 
\begin{itemize}
\item[]{\bf $C$ is chosen so that $r$ is minimal such that~(\ref{C_i}) holds.}
\end{itemize}
As explained in Remark~3.9 of \cite{Djo06}, the number $r$ is an invariant of $Th(\mcM)$.
The next lemma shows that the terminology ``the height is 1'' in the sense of \cite{Djo06} 
is equivalent with saying that ``the height is 1'' in the sense of Definition~\ref{definition of height}.

\begin{lem}\label{height and existence of a coordinatizing rank one set}
The following are equivalent:
\begin{itemize}
\item[(1)] $r = 1$.
\item[(2)] There is a $\es$-definable $D \subseteq M\meq$ in which only finitely many sorts are represented
such that $M \subseteq \acl(D)$ and $\su(d) = 1$ for every $d \in D$.
\end{itemize}
\end{lem}

\noindent
{\bf Proof.}
The implication from~(1) to~(2) is immediate, because if (1) holds, then $M \subseteq \acl(C_1)$
where $\su(c) = 1$ for every $c \in C_1$ by~(\ref{C_i}).

For the other direction, suppose that~(2) holds.
Then take $C_0 = \es$, $C_1 = D$ and $C = C_1 \cup M$.
Now it is straightforward to  verify that~(\ref{C_i}) holds for $r = 1$.
\hfill $\square$

\begin{lem}\label{r=1}
If $\mcM$ is primitive, then $r = 1$, so $M \subseteq \acl(C_1)$.
\end{lem}

\noindent
{\bf Proof.}
Suppose that $\mcM$ is primitive.
Towards a contradiction, suppose that $r > 1$.
Consider the following equivalence relation on $M$:
\[ x \sim y \ \Longleftrightarrow \ \acl(x) \cap C_{r-1} = \acl(y) \cap C_{r-1}. \]
By Fact~\ref{facts about simple omega-categorical structures}, this relation is $\es$-definable in $\mcM$.
By the assumption that $r > 1$ and that $r$ is minimal such that~(\ref{C_i}) holds,
there is $a \in M$ such that $a \notin \acl(C_{r-1})$.
Hence $a \notin \acl(\acl(a) \cap C_{r-1})$.
Then (by Fact~\ref{facts about simple omega-categorical structures})
there is $a' \in M$ such that $a' \neq a$ and
\[ \tp(a' / \acl(a) \cap C_{r-1}) = \tp(a / \acl(a) \cap C_{r-1}). \]
It follows that $\acl(a') \cap C_{r-1} = \acl(a) \cap C_{r-1}$, so $a \sim a'$.
By the assumption that $\mcM$ is primitive, we must have $b \sim a$ for all $b \in M$.
In other words,
\begin{equation}\label{all element have the same acl in C_r-1}
\text{ for all } b \in M, \ \acl(b) \cap C_{r-1} = \acl(a) \cap C_{r-1}.
\end{equation}
Let $A = \acl(a) \cap C_{r-1}$, so $A$ is finite; let $\bar{a}$ enumerate $A$.
Since, by~(\ref{C_i}) (and properties of dividing), $\su(b) \geq 1$ for every $b \in A$,
there is $\bar{b} \in C_{r-1}$ such that $\rng(\bar{b}) \neq A$ and
$\tp(\bar{b}) = \tp(\bar{a})$.
Hence there is $a'' \in M$ such that
\[ \acl(a'') \cap C_{r-1} = \rng(\bar{b}) \neq A = \acl(a) \cap C_{r-1}, \]
which contradicts~(\ref{all element have the same acl in C_r-1}).
Hence $r = 1$ and now~(\ref{C_i}) immediately gives that $M \subseteq \acl(C_1)$.
\hfill $\square$
\\

\noindent
Our aim is to prove Theorems~\ref{first main result} and~\ref{second main result} where it is assumed that
$\mcM$ is primitive or has weight 1.
It follows from the definition of having height 1 (Definition~\ref{definition of height}) and 
Lemmas~\ref{height and existence of a coordinatizing rank one set}
and~\ref{r=1} that in either case we have $r=1$. Therefore

\begin{itemize}
\item[] {\bf we assume for the rest of this section that $r = 1$, so $M \subseteq \acl(C_1)$.}
\end{itemize}

\noindent
By~(\ref{C_i}),  for every $c \in C_1$, $\su(c) = 1$.
As dependence is trivial it follows that if $c \in C_1$, $A \subseteq C_1$ and $c \in \acl(A)$, then $c \in \acl(a)$ for some $a \in A$.
Now consider the equivalence relation `$\acl(x) \cap C_1 = \acl(y) \cap C_1$' on $C_1$.
By Fact~\ref{facts about simple omega-categorical structures} this relation is $\es$-definable and there is $t < \omega$ such
that every equivalence class has at most $t$ elements.
Since $\mcM\meq$ has elimination of imaginaries (see \cite{Hod, She})
it follows  that each equivalence class corresponds
to an element in $\mcM\meq$ in the following sense:
There is a $\es$-definable set $C'_1 \subseteq M\meq$ in which only finitely many sorts are represented 
and a $\es$-definable surjective function $f : C_1 \to C'_1$ (meaning that the graph of $f$
is a $\es$-definable relation) such that if $c, c' \in C_1$ then $f(c) = f(c')$ if and only if $\acl(c) \cap C_1 = \acl(c') \cap C_1$.
It follows that for every $c \in C_1$, $f(c) \in \dcl(c)$ and $c \in \acl(f(c))$.
This implies that $M \subseteq \acl(C'_1)$ and if $c \in C'_1$, then $\su(c) = 1$ and $\acl(c) \cap C'_1 = \{c\}$. 
Since dependence is trivial, it follows that $\acl(A) \cap C'_1 = A$ for every $A \subseteq C'_1$.
In order not to switch from the notation `$C_1$' to the notation `$C'_1$', we may (by the above argument), without loss of generality, assume that
\begin{equation}\label{C_1 is a trivial geometry}
\text{ for every } A \subseteq C_1, \ \acl(A) \cap C_1 = A.
\end{equation}

\begin{defin}\label{definition of crd}{\rm
For every $a \in M\meq$, let $\crd(a) = \acl(a) \cap C_1$. We call `$\crd(a)$' the (set of) {\em coordinates} of $a$.
(Our definition of $\crd$ corresponds to the notation $\crd_1$ in \cite{Djo06}.)
}\end{defin}

\noindent
From the definition we obviously have  $\crd(a) \subseteq \acl(a)$ for every $a \in M$.

\begin{lem}\label{every element is in dcl of its coordinates}
(i) For every $a \in M$, $a \in \acl(\crd(a))$ and hence $\acl(a) = \acl(\crd(a))$.\\
(ii) For all $a, a' \in M$, $\acl_\mcM(a) = \acl_\mcM(a')$ if and only if $\crd(a) = \crd(a')$.\\
(iii) If the equivalence relation $\acl_\mcM(x) = \acl_\mcM(y)$ has only singleton classes, 
then $a \in \dcl(\crd(a))$ for every $a \in M$.
\end{lem}

\noindent
{\bf Proof.}
Part~(i) is a direct consequence of  Lemma~5.1 in \cite{Djo06}, but can easily be proved directly as follows.
Let $a \in M$, so $a \in \acl(C_1)$.
Let $B = \acl(a) \cap C_1$.
Suppose that $a \notin \acl(B)$. 
Then $a \underset{B}{\nind} C_1$.
As dependence is trivial there is $c \in C_1 \setminus B$ such that 
$a \underset{B}{\nind} c$ and hence $aB \nind c$
Since $B \subseteq \acl(a)$ we get $a \nind c$ which (as $\su(c) = 1$) implies that
$c \in \acl(a)$. Hence $c \in B$ which contradicts the choice of $c$.

Part~(ii) follows immediately from part~(i) and the definition of $\crd$.
Hence it remains to prove~(iii).

Suppose that $a, a' \in M$ and $a' \in \acl(\crd(a))$.
By the primitivity of $\mcM$ we have $\tp(a) = \tp(a')$, hence $|\crd(a)| = \crd(a')|$ and therefore 
$\crd(a) = \crd(a')$. Hence $a$ and $a'$ belong to the same equivalence class
of $\acl_\mcM(x) = \acl_\mcM(y)$ and thus $a = a'$.
\hfill $\square$

\begin{rem}\label{remark on lemma about properties of crd}{\rm
Let `$\approx$' be the ($\es$-definable) equivalence relation `$\acl_\mcM(x) = \acl_\mcM(y)$' and 
$\mcM / \approxc$ the canonically embedded structure of $\mcM\meq$ with universe $M / \approxc$.
For every $a \in M$ we have $\acl(a) = \acl([a]_\approx)$.
Therefore $\crd(a) = \crd([a]_\approx)$ for every $a \in M$ and Lemma~\ref{every element is in dcl of its coordinates}
holds if `$M$'  and '$\mcM$' are replaced with `$M / \approxc$' and `$\mcM / \approxc$', respectively.
Since the equivalence relation `$\acl_{\mcM / \approx}(x) = \acl_{\mcM / \approx}(y)$ is trivial
it follows that part~(iii) simplifies to the statement: 
for every $a \in M / \approxc$, $a \in \dcl(\crd(a))$.
}\end{rem}

\begin{lem}\label{crd is unique up to isomorphism}
If $a, b \in M$ and $\tp(a) = \tp(b)$, then $\crd(a)$ and $\crd(b)$ can be ordered as $\bar{a}$ and $\bar{b}$, respectively,
so that $\tp(\bar{a}) = \tp(\bar{b})$.
\end{lem}

\noindent
{\bf Proof.}
Suppose that $a, b \in M$ and $\tp(a) = \tp(b)$.
Let $\bar{a}$ be an ordering of $\crd(a)$ and let $f$ be an automorphism of $\mcM\meq$
such that $f(a) = b$. Then $f(\bar{a})$ is an ordering of $\crd(b)$ such that
$\tp(\bar{a}) = \tp(f(\bar{a}))$.
\hfill $\square$
\\

\noindent
By considering a $\es$-definable subset of $C_1$ if necessary, we may, in addition to previous assumptions, assume that
\begin{equation}\label{every member of C_1 is a coordinate of M}
\text{ for every $c \in C_1$, there is $a \in M$ such that $c \in \crd(a)$.}
\end{equation}

\noindent
Let $\mcC_1$ be the canonically embedded structure in $\mcM\meq$ with universe  $C_1$.
By Fact~\ref{rank 1 sets are reducts of binary random structures}, 
$\mcC_1$ is a reduct of a binary random structure $\mcR$, so in particular $R = C_1$, where $R$ is the universe of $\mcR$.
Hence $\crd(a) \subseteq R$ for every $a \in M$.

\begin{prop}\label{M is definable in Req}
Let `$\approx$' be the $\es$-definable equivalence relation on $M$ defined by $x \approx y$ if and only if $\acl_\mcM(x) = \acl_\mcM(y)$.
Then $\mcM / \approxc$ is definable in $\mcR\meq$  and we can use a formula $\chi(x)$ and bijection $f : \chi(\mcR\meq) \to M / \approxc$
as in Definition~\ref{definition of strong interpretability} with the properties that 
\begin{enumerate}
\item for all $c \in \chi(\mcR\meq)$,
$\acl_{\mcR\meq}(c) \cap R = \crd(f(c))$ and $\acl_{\mcR\meq}(c) = \acl_{\mcR\meq}(\crd(f(c)))$, 

\item for all $c, c' \in \chi(\mcR\meq)$, $c = c'$ if and only if $\crd(f(c)) = \crd(f(c'))$, and

\item $c \in \dcl_{\mcR\meq}(\acl_{\mcR\meq}(c) \cap R)$ for all $c \in \chi(\mcR\meq)$.
\end{enumerate}
\end{prop}

\noindent
{\bf Proof.}
In order to simplify notation and make the argument more clear we prove the proposition under the extra assumption
that $\approx$ is  trivial, in other words, that it has only singleton classes. The general case is a straightforward modification of this special case,
where we use Remark~\ref{remark on lemma about properties of crd} instead of Lemma~\ref{every element is in dcl of its coordinates}.
The assumption that $\approx$ is trivial implies that $\mcM / \approxc$ is definable in $\mcM$ via the map taking every $a \in M$
to $[a]_\approx$, where $[a]_\approx$ is the (singleton) $\approx$-class to which $a$ belongs.
Hence it suffices to show that $\mcM$ is definable in $\mcR\meq$.

Observe that since $\mcC_1$ is a reduct of $\mcR$ we can (and will) view $(\mcC_1)\meq$ as a reduct of $\mcR\meq$.
Also, since $\mcC_1$ is canonically embedded in $\mcM\meq$ we can (and will) assume that 
$(C_1)\meq \subseteq (M\meq)\meq$. 
Moreover, as $\mcC_1$ is canonically embedded in $\mcM\meq$,
for all $\bar{c}, \bar{c}' \in (C_1)\meq$, $\tp_{(\mcC_1)\meq}(\bar{c}) = \tp_{(\mcC_1)\meq}(\bar{c}')$
if and only if $\tp_{(\mcM\meq)\meq}(\bar{c}) = \tp_{(\mcM\meq)\meq}(\bar{c}')$.

Let $p_1, \ldots, p_s$ enumerate (without repetition) $S^\mcM_1(\es)$.
For each $p_i$ choose a realization $a_i \in M$ of $p_i$ and enumerate (without repetition) $\crd(a_i)$ as $\bar{b}_i$.
Then let $q_i = \tp(\bar{b}_i)$.
Since $\mcM$ is $\omega$-categorical and since $\mcC_1$ is canonically embedded in $\mcM\meq$ and
$\mcC_1$ is a reduct of $\mcR$ it follows that (for each $i$) the set $q_i(\mcM\meq)$ is $\es$-definable in $\mcM\meq$,
$\mcC_1$ and in $\mcR$.

For each $i$, define an equivalence relation on $q_i(\mcM\meq)$ as follows: 
$\bar{x} \sim_i \bar{y}$ if and only if $\rng(\bar{x}) = \rng(\bar{y})$.
(The relation $\sim_i$ can be extended to all $|\bar{b}_i|$-tuples of elements from $M\meq$, $R$ or $C_1$,
by saying that all $|\bar{b}_i|$-tuples outside of $q_i(\mcM\meq)$ belong to the same class.)
Note that $\sim_i$ is a $\es$-definable relation in $\mcR$ as well as in $\mcC_1$ and in $\mcM\meq$.
Hence, for every $i$, the set of $\sim_i$-classes of tuples in $q_i(\mcM\meq)$ is a $\es$-definable subset of $\mcR\meq$, of $(\mcC_1)\meq$ and of
$(\mcM\meq)\meq$.

For every $i$ and $\bar{b}$ realizing $q_i$ let $[\bar{b}]_i$ be its $\sim_i$-class.
Then let
\[X \ = \ \{ c : c = [\bar{b}]_i \text{ for some $i$ and $\bar{b}$} \} \]
so $X$ is a subset of $R\meq$, $(C_1)\meq$ and of $(M\meq)\meq$.
Moreover, $X$ is $\es$-definable in $\mcR\meq$ by some formula $\chi(x)$. 

Now we define a bijection $g : M \to X$ such that if $f = g^{-1}$ then $f$ has the required properties.
For every  $a \in M$ define $g(a)$ as follows: 
let $i$ be such that $\mcM \models p_i(a)$ and (using Lemma~\ref{crd is unique up to isomorphism})
let $\bar{b}$ enumerate $\crd(a)$ in such a way that
$\bar{b}$ realizes $q_i$ and let $g(a) = [\bar{b}]_i$.
The surjectivity of $g$ follows from the $\omega$-homogeneity of $\mcM\meq$ and the definition of $X$.
Observe that if $a \in M$ realizes $p_i$ and $\bar{b}$ enumerates $\crd(a)$ in such a way that it realizes $q_i$
then (by Lemma~\ref{every element is in dcl of its coordinates})
$\acl(\bar{b}) = \acl(a)$ and $[\bar{b}]_i \in \dcl_{(\mcM\meq)\meq}(a)$. 
Since each $\sim_i$-class is finite we get 
$\acl_{(\mcM\meq)\meq}([\bar{b}]_i) = \acl_{(\mcM\meq)\meq}(a)$.
As we assume that the equivalence relation $\acl_\mcM(x) = \acl_\mcM(y)$ has only singleton classes 
(and hence the same holds for $\acl_{(\mcM\meq)\meq}(x) = \acl_{(\mcM\meq)\meq}(y)$ restricted to $M$) it follows that 
$g : M \to X$ is bijective and $a \in \dcl_{(\mcM\meq)\meq}(g(a))$ for all $a \in M$. 
(In the general case, note that the equivalence relation $\acl_{\mcM / \approx}(x) = \acl_{\mcM / \approx}(y)$
is trivial, by the definition of $\approx$.)

Let $0 < n < \omega$ and $a_1, \ldots, a_n, a'_1, \ldots, a'_n \in M$.
By the observations already made (following from the fact that $\mcC_1$ is a reduct of $\mcR$ and $\mcC_1$ is
canonically embedded in $\mcM\meq$) and in particular since
 $a_i \in \dcl_{(\mcM\meq)\meq}(g(a_i))$ for all $i$ and similarly for $a'_i$
we get:
\begin{align*}
\tp_{\mcR\meq}(g(a_1), \ldots, g(a_n)) \ &= \ \tp_{\mcR\meq}(g(a'_1), \ldots, g(a'_n)) \ \ \Longrightarrow \\
\tp_{(\mcC_1)\meq}(g(a_1), \ldots, g(a_n)) \ &= \ \tp_{(\mcC_1)\meq}(g(a'_1), \ldots, g(a'_n)) \ \ \Longrightarrow \\
\tp_{(\mcM\meq)\meq}(g(a_1), \ldots, g(a_n)) \ &= \ \tp_{(\mcM\meq)\meq}(g(a'_1), \ldots, g(a'_n)) \ \ \Longrightarrow \\
\tp_{(\mcM\meq)\meq}(a_1, \ldots, a_n) \ &= \ \tp_{(\mcM\meq)\meq}(a'_1, \ldots, a'_n) \ \ \Longrightarrow \\
\tp_\mcM(a_1, \ldots, a_n) \ &= \ \tp_\mcM(a'_1, \ldots, a'_n).
\end{align*}
As $\mcR$ is homogeneous, every type over $\es$ with respect to $\mcR\meq$ which is realized by elements from $X$ is isolated.
Moreover, for each $0 < n < \omega$, only finitely many types from $S^{\mcR\meq}_n(\es)$ are realized by $n$-tuples from $X^n$.
It follows that for every $0 < n < \omega$ and $p \in S^\mcM_n(\es)$ there is a formula $\varphi_p(x_1, \ldots, x_n)$ without parameters
in the language of $\mcR\meq$ such that for all $a_1, \ldots, a_n \in M$,
\[ \mcM \models p(a_1, \ldots, a_n) \ \ \Longleftrightarrow \ \ \mcR\meq \models \varphi_p(g(a_1), \ldots, g(a_n)). \]
This implies that $\mcM$ is definable in $\mcR\meq$ via the map $f = g^{-1}$ as in Definition~\ref{definition of strong interpretability}.

It remains to verify that $f$ has the other properties stated in the proposition.
Let $c \in \chi(\mcR\meq)$, so $c = [\bar{b}]_i$ for some $i$ and some $\bar{b} \in q_i(\mcM\meq)$.
From the definition of $f$ it follows that $\rng(\bar{b}) = \crd(f(c))$, so $\acl_{\mcR\meq}(\bar{b}) = \acl_{\mcR\meq}(\crd(f(c)))$.
As $c = [\bar{b}]_i$ is a finite equivalence class we get
$\acl_{\mcR\meq}(c) = \acl_{\mcR\meq}(\bar{b}) = \acl_{\mcR\meq}(\crd(f(c)))$.
In particular,
$\acl_{\mcR\meq}(c)  \cap R = \acl_{\mcR\meq}(\bar{b}) \cap R$, and since $\bar{b} \in R$ where
$\mcR$ is a binary random structure we get $\acl_{\mcR\meq}(\bar{b}) \cap R = \rng(\bar{b}) = \crd(f(c))$.
Hence $\acl_{\mcR\meq}(c) \cap R = \crd(f(c))$. This proves~(1).
Since $f$ is bijective and the relation $\acl_\mcM(x) = \acl_\mcM(y)$ has only singleton classes (by assumption) it follows 
(using Lemma~\ref{every element is in dcl of its coordinates}~(ii)) that
for all $c, c' \in \chi(\mcR\meq)$, $c = c'$ if and only if $\crd(f(c)) = \crd(f(c'))$, so~(2) is proved.
For~(3), suppose that $c, c' \in \chi(\mcR\meq)$ are distinct and that $c' \in \acl_{\mcR\meq}(\acl_{\mcR\meq}(c) \cap R) =
\acl_{\mcR\meq}(\crd(f(c)) = \acl_{\mcR\meq}(c)$. Then $\crd(f(c')) = \acl_{\mcR\meq}(c') \cap R \subseteq 
\acl_{\mcR\meq}(c) \cap R = \crd(f(c))$, so by~(2), $|\acl_{\mcR\meq}(c') \cap R| < |\acl_{\mcR\meq}(c) \cap R|$. 
This implies that $\tp_{\mcR\meq}(c') \neq \tp_{\mcR\meq}(c)$ whenever 
$c, c' \in \chi(\mcR\meq)$ are distinct and $c' \in \acl_{\mcR\meq}(\acl_{\mcR\meq}(c) \cap R)$.
It follows that $c \in \dcl_{\mcR\meq}(\acl_{\mcR\meq}(c) \cap R)$.
\hfill $\square$
\\

\noindent
{\em Observe that Proposition~\ref{M is definable in Req} proves part~(i) of Theorem~\ref{second main result}.}
We continue by proving part~(ii) of Theorem~\ref{second main result}.
Therefore,
\begin{itemize}
\item[] {\bf we assume for the rest of the section that the equivalence relation $\acl_\mcM(x) = \acl_\mcM(y)$ is trivial.}
\end{itemize}
Then there is an obvious bijection $h : M / \approxc \ \to M$ such that for all 
$\bar{a}, \bar{b} \in M / \approxc$, $\tp(\bar{a}) = \tp(\bar{b})$ if and only if
$\tp(h(\bar{a})) = \tp(h(\bar{b}))$.
Therefore Proposition~\ref{M is definable in Req} allows us to identify $M$ with $\chi(\mcR\meq)$ via the bijection $h \circ f$ where 
$\chi$ and $f$ are as in that proposition.
It follows that for every $0 < n < \omega$ and $D \subseteq M^n$, if $D$ is $\es$-definable in $\mcM$, then
it is $\es$-definable in $\mcR\meq$. 
So 
\begin{itemize}
\item[] {\bf for the rest of this section we assume that $\mcM$ is a reduct of the canonically embedded structure of $\mcR\meq$ 
with universe $M = \chi(\mcR\meq)$.}
\end{itemize}
Note that the identification $h(f(c)) = c$
for all $c \in \chi(\mcR\meq)$ together with~(1) of Proposition~\ref{M is definable in Req} implies that
\[\text{ for all $a \in M$}, \ \crd(a) \ = \ \acl_{\mcR\meq}(a) \cap R.\]
Part~(3) of Proposition~\ref{M is definable in Req} and the new assumptions imply that
\begin{equation}\label{a is in definable closure of acl(a) cut with R}
\text{for every $a \in M$, $a \in \dcl_{\mcR\meq}(\crd(a))$.}
\end{equation}

\begin{lem}\label{crd(a) is rigid}
Let $a \in M$ and $\crd(a)  = \{b_1, \ldots, b_m\}$, where the elements are enumerated without repetition.
Then for every nontrivial permutation $\pi$ of $\{1, \ldots, m\}$,
$\tp_\mcR(b_1, \ldots, b_m) \neq \tp_\mcR(b_{\pi(1)}, \ldots, b_{\pi(m)})$.
\end{lem}

\noindent
The proof of Lemma~\ref{crd(a) is rigid} is given in Section~\ref{Proof of the main technical result},
but now we derive a corollary of it.

\begin{cor}\label{crd is in dcl of a}
For every  $a \in M$, $\crd(a) \subseteq \dcl_{\mcR\meq}(a)$.
\end{cor}

\noindent
{\bf Proof.}
Let $a \in M$ and $\crd(a) = \{b_1, \ldots, b_m\}$.
If $b_i \notin \dcl_{\mcM\meq}(a)$ for some $i$, then 
$\tp_{\mcR\meq}(a, b_i) = \tp_{\mcR\meq}(a, b_j)$ for some $j \neq i$, from which it follows
(using that $\mcR\meq$ is $\omega$-homogeneous) that
there is a nontrivial permutation $\pi$ of $\{1, \ldots, m\}$ such that 
\[ \tp_{\mcR\meq}(a, b_1, \ldots, b_m) \ = \  \tp_{\mcR\meq}(a, b_{\pi(1)}, \ldots, b_{\pi(m)}). \]
Then $\tp_\mcR(b_1, \ldots, b_m) = \tp_\mcR(b_{\pi(1)}, \ldots, b_{\pi(m)})$ which contradicts
Lemma~\ref{crd(a) is rigid}.
\hfill $\square$
\\

\noindent
Now we are ready to prove the remaining parts of the main results (and when stating them we repeat the assumptions made
in the beginning of this section).

\bigskip

\noindent
{\bf Theorem~\ref{second main result}}{\em \ 
Suppose that $\mcM$ is a structure which is countable, binary, homogeneous, simple and 1-based with height 1.
Let $\approx$ denote the equivalence relation `$\acl_\mcM(x) = \acl_\mcM(y)$'.\\
(i) Then there is a binary random structure $\mcR$ such that $\mcM / \approxc$ is definable in $\mcR\meq$.\\
(ii) If, in addition, the equivalence relation `$\approx$' is trivial, 
then $\mcM$ is strongly interpretable in a binary random structure.
}

\medskip
\noindent
{\bf Proof.}
As mentioned above, part~(i) follows from Proposition~\ref{M is definable in Req}.
So it remains to prove~(ii) and for this we adopt the assumption that the equivalence relation $\approx$ is trivial, as well as
all other assumptions that have been made earlier in this section.
The proof is similar to the proof of Proposition~\ref{M is definable in Req}. 
The essential difference is that the assumption that $\approx$ is trivial together with Lemma~\ref{crd(a) is rigid} allows
us to reach a stronger conclusion than in Proposition~\ref{M is definable in Req}.

Let $p_1, \ldots, p_s$ enumerate $S^\mcM_1(\es)$.
For each $1 \leq i \leq s$ choose a realization $a_i \in M$ of $p_i$ and then choose an ordering
$\bar{b}_i$ of $\crd(a_i)$.
For each $i$, let $\theta_i$ isolate $\tp_{\mcR\meq}(a_i, \bar{b}_i)$ and let
$\varphi_i$ isolate $\tp_\mcR(\bar{b}_i)$.
Let $X = \varphi_1(\mcR) \cup \ldots, \cup \varphi_s(\mcR)$.

By~(\ref{a is in definable closure of acl(a) cut with R}),
for every $\bar{b} \in X$ there is a unique $a \in M$ such that 
$\mcR\meq \models \theta_i(a, \bar{b})$ for some $i$. 
For every $\bar{b} \in X$ we let $f(\bar{b}) = a$ for the unique $a \in M$ such that $\mcR\meq \models \theta_i(a, \bar{b})$ for some $i$.
Since $\crd(a)$ exists as a subset of $C_1 = R$ for every $a \in M$ it follows that $f : X \to M$ is surjective.
Now we claim that $f$ is injective. 
For if $\bar{b}, \bar{b}' \in X$ and $f(\bar{b}) = f(\bar{b}') = a \in M$,
then $\mcR\meq \models \theta_i(a, \bar{b}) \wedge \theta_j(a, \bar{b}')$ for some $i, j$.
Then $\tp_{\mcR\meq}(a, \bar{b}) = \tp_{\mcR\meq}(a_i, \bar{b}_i)$ and
$\tp_{\mcR\meq}(a, \bar{b}') = \tp_{\mcR\meq}(a_j, \bar{b}_j)$ from which it follows
that $\tp_{\mcR\meq}(a_i) = \tp_{\mcR\meq}(a_j)$ and $\rng(\bar{b}) = \crd(a)  = \rng(\bar{b}')$. 
From the construction it now follows that $i = j$ and hence
$\tp_{\mcR\meq}(a, \bar{b}) = \tp_{\mcR\meq}(a, \bar{b}')$.
By Corollary~\ref{crd is in dcl of a}, $\rng(\bar{b}) \subseteq \dcl_{\mcR\meq}(a)$ and $\rng(\bar{b}') \subseteq \dcl_{\mcR\meq}(a)$
so we must have $\bar{b} = \bar{b}'$.

Let $0 < n < \omega$, $1 \leq j \leq s$ and $\bar{b}_1, \ldots, \bar{b}_n, \bar{b}'_1, \ldots, \bar{b}'_n \in \varphi_j(\mcR)$.
From the definition of $f$, its graph is a $\es$-definable relation, which implies 
(since $f$ is bijective between two $\es$-definable subsets of $R\meq$)
that 
$\rng(\bar{b}_i) \subseteq \dcl_{\mcR\meq}(f(\bar{b}_i))$ and $f(\bar{b}_i) \in \dcl_{\mcR\meq}(\bar{b}_i)$ for each $i$, and similarly for
each $\bar{b}'_i$. This implies that
\begin{align*}
\tp_\mcR(\bar{b}_1, \ldots, \bar{b}_n) \ &= \ \tp_\mcR(\bar{b}'_1, \ldots, \bar{b}'_n) \ \Longleftrightarrow \\
\tp_{\mcR\meq}(\bar{b}_1, \ldots, \bar{b}_n) \ &= \ \tp_{\mcR\meq}(\bar{b}'_1, \ldots, \bar{b}'_n) \ \Longleftrightarrow \\
\tp_{\mcR\meq}(f(\bar{b}_1), \ldots, f(\bar{b}_n)) \ &= \ \tp_{\mcR\meq}(f(\bar{b}'_1), \ldots, f(\bar{b}'_n)) \ \Longrightarrow \\
\tp_\mcM(f(\bar{b}_1), \ldots, f(\bar{b}_n)) \ &= \ \tp_\mcM(f(\bar{b}'_1), \ldots, f(\bar{b}'_n)).
\end{align*}
Since $\mcR$ and $\mcM$ are homogeneous and hence $\omega$-categorical,
there is 
for every $0 < n < \omega$ and formula $\xi(x_1, \ldots, x_n)$ in the language of $\mcM$
a formula $\xi'(\bar{x}_1, \ldots, \bar{x}_n)$ in the language of $\mcR$ such that for
all $\bar{b}_1, \ldots, \bar{b}_n \in \varphi_i(\mcR)$,
$\mcM \ \models \xi(f(\bar{b}_1), \ldots, f(\bar{b}_n))$ if and only if
$\mcR \models \xi'(\bar{b}_1, \ldots, \bar{b}_n)$.
By Lemma~\ref{equivalent condition of being strongly interpretable}, 
$\mcM$ is strongly interpretable in $\mcR$.
\hfill $\square$

\bigskip

\noindent
{\bf Theorem~\ref{first main result}}{\em \ 
Suppose that $\mcM$ is a structure which is countable, binary, homogeneous, primitive, simple and 1-based.
Then $\mcM$ is strongly interpretable in a  binary random structure.
}

\medskip
\noindent
{\bf Proof.}
Suppose that $\mcM$ satisfies the assumptions of the Theorem~\ref{first main result}. The primitivity of $\mcM$ implies that
the equivalence relation $\acl_\mcM(x) = \acl_\mcM(y)$ is trivial, because it is $\es$-definable in $\mcM$.
Therefore Theorem~\ref{first main result} is a direct consequence of 
Lemma~\ref{r=1} and 
Theorem~\ref{second main result}.
\hfill $\square$

\begin{examps}\label{remark on construction of line structures}{\rm
One can construct structures that satisfy the assumptions of Theorem~\ref{second main result}
roughly as follows. Take a binary random structure $\mcR$.
Choose some rigid nonisomorphic substructures $\mcL_1, \ldots, \mcL_k$ of $\mcR$ (where {\em rigid} means that there is no
nontrivial automorphism) together with a fixed enumeration of $L_i$ for each $i$.
Let $M$ be the set of substructures of $\mcR$ which are isomorphic to some $\mcL_i$.
If $\mcL \cong \mcL_i$ then we think of $L$ as being enumerated so that this enumeration
and the enumeration of $L_i$ induces an isomorphism from $\mcL$ to $\mcL_i$.
For any $\mcA, \mcA', \mcB, \mcB' \in M$, let $(\mcA, \mcA')$ and $(\mcB, \mcB')$
have the same type in $\mcM$
if and only if $\bar{a}\bar{a}'$ and $\bar{b}\bar{b}'$ have the same type in $\mcR$ where 
$\bar{a}, \bar{a}', \bar{b}, \bar{b}'$ are the enumerations of the elements in the respective structure.
By Theorem~\ref{first main result} and 
Lemma~\ref{crd(a) is rigid} 
it follows that if $\mcN$ is countable, binary, homogeneous, primitive,
simple and 1-based, then $\mcN$ is a proper reduct of such a structure $\mcM$ (with $k = 1$).
}
\end{examps}

\section{Proof of Lemma~\ref{crd(a) is rigid}}\label{Proof of the main technical result}

\noindent
In this section we prove Lemma~\ref{crd(a) is rigid}, so all assumptions in Section~\ref{Binary homogeneous 1-based structures}
up to the Lemma~\ref{crd(a) is rigid} apply in this section, including the conventions of Notation~\ref{notation for acl etc}.
In particular we recall that $\mcR$ is a binary random structure and $\mcC_1$ is the canonically embedded structure of $\mcM\meq$
with universe $C_1$.

The intuition behind the proof of Lemma~\ref{crd(a) is rigid} comes from the argument that the line graph over
a complete graph with infinite countable vertex set is not homogeneous. 
Somewhat more precisely: Suppose that $\mcK$ is a complete graph with infinite countable vertex set $K$.
The line graph over $\mcK$ has as its vertex set the set of all 2-subsets of $K$ and two 2-subsets are adjacent (in the line graph)
if and only if they intersect in exactly one point. Then one can choose distinct $a, b, c \in K$
and distinct $d, e, e', e'' \in K$ and it is easy to see that the 3-tuples
$(\{a, b\}, \{b, c\}, \{c, a\})$ and $(\{d, e\}, \{d, e'\}, \{d, e''\})$ satisfy the same quantifier free 
formulas in the line graph, but there is a formula which is satisfied by one of the tuples but not the other.

\medskip

\noindent
{\bf Lemma~\ref{crd(a) is rigid}}{\em \
Let $a \in M$ and $\crd(a)  = \{b_1, \ldots, b_m\}$, where the elements are enumerated without repetition.
Then for every nontrivial permutation $\pi$ of $\{1, \ldots, m\}$,
$\tp_\mcR(b_1, \ldots, b_m) \neq \tp_\mcR(b_{\pi(1)}, \ldots, b_{\pi(m)})$.
}

\medskip

\noindent
{\bf Proof.}
As usual we use the facts from Section~\ref{Preliminaries} without further reference.
The lemma is trivial if $m = 1$, so we assume that $m \geq 2$.
Let $a \in M$ and let 
\[ \crd(a)  = \{b_1, \ldots, b_m\}. \]
Towards a contradiction, suppose that there is a nontrivial permutation $\pi$ of $\{1, \ldots, m\}$
such that 
\begin{equation}\label{the permuted sequence has the same type}
\tp_\mcR(b_1, \ldots, b_m) = \tp_\mcR(b_{\pi(1)}, \ldots, b_{\pi(m)}).
\end{equation}
To simplify notation and witout loss of generality,
\begin{equation}\label{pi(1) equals 2}
\text{ we assume that $\pi(1) = 2$. }
\end{equation}
From~(\ref{the permuted sequence has the same type}) and~(\ref{pi(1) equals 2}) we get
\[ \tp_\mcR(b_1) = \tp_\mcR(b_2). \]
As $\mcR$ is a binary random structure we can argue as in 
Remark~\ref{remarks about binary random structures}
(with $n = 3$, $p_{1, 2}(x_1, x_2) = \tp_\mcR(b_1, b_2)$, 
$p_{1, 3}(x_1, x_3) = \tp_\mcR(b_1, b_2)$ and $p_{2, 3}(x_2, x_3) = \tp_\mcR(b_1, b_2)$)
and find $b'_2 \in R$ such that 
\[ \tp_\mcR(b_1, b'_2) \ = \ \tp_\mcR(b_1, b_2) \ = \ \tp_\mcR(b_2, b'_2). \]
If $m > 2$, then, by using that $\mcR$ is a random structure again, we find 
{\em distinct} elements $b'_3, \ldots, b'_m, b''_3, \ldots, b''_m \in R$ such that if we let
\begin{equation}\label{identifications of elements}
b'_1 = b_1, \ b''_1 = b_2, \ b''_2 = b'_2,
\end{equation}
then
\begin{align}\label{the other sequences have the same type}
\tp_\mcR(b_1, b_2, \ldots, b_m) \ = \ 
\tp_\mcR(b'_1, b'_2, \ldots, b'_m) \ = \ 
\tp_\mcR(b''_1, b''_2, \ldots, b''_m). 
\end{align}

\noindent
From~(\ref{the permuted sequence has the same type}) and~(\ref{the other sequences have the same type})
we get
\begin{align}\label{the same type when the other sequences are permuted}
\tp_\mcR(b_1, b_2, \ldots, b_m) \ &= \ \tp_\mcR(b_{\pi(1)}, b_{\pi(2)}, \ldots, b_{\pi(m)}),  \\
\tp_\mcR(b'_1, b'_2, \ldots, b'_m) \ &= \ \tp_\mcR(b'_{\pi(1)}, b'_{\pi(2)}, \ldots, b'_{\pi(m)}) \quad \text{ and} \nonumber \\
\tp_\mcR(b''_1, b''_2, \ldots, b''_m) \ &= \ \tp_\mcR(b''_{\pi(1)}, b''_{\pi(2)}, \ldots, b''_{\pi(m)}), \nonumber
\end{align}
Since $\tp_\mcR(b_1) = \tp_\mcR(b_2) = \tp_\mcR(b'_2)$ and $\mcR$ is a binary random structure it follows 
from~(\ref{the same type when the other sequences are permuted}) that there are
$c_2, \ldots, c_m \in R \setminus \{b_1, \ldots, b_m, b'_1, \ldots, b'_m\}$ such that if $c_1 = b_1 (= b'_1)$,
then
\begin{align}\label{equal types of some long sequences}
\tp_\mcR(b_1, \ldots, b_m,  c_1, \ldots, c_m)  \ &= \ \tp_\mcR(b_{\pi(1)}, \ldots, b_{\pi(m)},  b''_1, \ldots, b''_m)
\quad \text{ and} \\
\tp_\mcR(b'_1, \ldots, b'_m, c_1, \ldots, c_m) \ &= \ \tp_\mcR(b'_{\pi(1)}, \ldots, b'_{\pi(m)}, b''_{\pi(1)}, \ldots, b''_{\pi(m)}).
\nonumber
\end{align}

\noindent
Since $\mcC_1$ is a reduct of $\mcR$ we can replace `$\tp_\mcR$' with `$\tp_{\mcC_1}$' everywhere 
 in~(\ref{the other sequences have the same type}),~(\ref{the same type when the other sequences are permuted})
and~(\ref{equal types of some long sequences}).
Moreover, as $\mcC_1$ is a canonically embedded structure of $\mcM\meq$, we can replace 
`$\tp_{\mcC_1}$' with `$\tp$' (which abbreviates `$\tp_{\mcM\meq}$'), so altogether we get
\begin{align}\label{type identifications in Meq}
&\tp(b_1, b_2, \ldots, b_m) \ = \ 
\tp(b'_1, b'_2, \ldots, b'_m) \ = \ 
\tp(b''_1, b''_2, \ldots, b''_m), \\
&\tp(b_1, b_2, \ldots, b_m) \ = \ \tp(b_{\pi(1)}, b_{\pi(2)}, \ldots, b_{\pi(m)}), \nonumber \\
&\tp(b'_1, b'_2, \ldots, b'_m) \ = \ \tp(b'_{\pi(1)}, b'_{\pi(2)}, \ldots, b'_{\pi(m)}), \nonumber \\
&\tp(b''_1, b''_2, \ldots, b''_m) \ = \ \tp(b''_{\pi(1)}, b''_{\pi(2)}, \ldots, b''_{\pi(m)}), \nonumber \\
&\tp(b_1, \ldots, b_m,  c_1, \ldots, c_m)  \ = \ \tp(b_{\pi(1)}, \ldots, b_{\pi(m)},  b''_1, \ldots, b''_m)
\quad \text{ and} \nonumber \\
&\tp(b'_1, \ldots, b'_m, c_1, \ldots, c_m) \ = \ \tp(b'_{\pi(1)}, \ldots, b'_{\pi(m)}, b''_{\pi(1)}, \ldots, b''_{\pi(m)}). \nonumber
\end{align}

\noindent
In particular,~(\ref{type identifications in Meq}) implies that 
$\tp(c_1, \ldots, c_m) = \tp(b''_{\pi(1)}, \ldots, b''_{\pi(m)}) = \tp(b''_1,  \ldots, b''_m) = \tp(b'_1, \ldots, b'_m) = \tp(b_1, \ldots, b_m)$.
So, using the $\omega$-homogeneity of $\mcM\meq$, there is $a^* \in M$ such that 
\[ \tp(a^*, c_1, \ldots, c_m) \ = \ \tp(a, b_1, \ldots, b_m) \]
and hence $\crd(a^*) = \{c_1, \ldots, c_m\}$.
For the same reason there is $a' \in M$ such that
\[ \tp(a', b'_1, \ldots, b'_m) \ = \ \tp(a, b_1, \ldots, b_m) \]
and hence $\crd(a') = \{b'_1, \ldots, b'_m\}$.
By the $\omega$-homogeneity of $\mcM\meq$ again
and~(\ref{type identifications in Meq}), there are $a_0, a'', a'_0, a''_0 \in M$ such that
\begin{align}\label{the elements corresponding to a}
\tp(a, b_1, \ldots, b_m, a^*, c_1, \ldots, c_m) \ &= \ \tp(a_0, b_{\pi(1)}, \ldots, b_{\pi(m)}, a'', b''_1, \ldots, b''_m) \ \ \text{and} \\
\tp(a', b'_1, \ldots, b'_m, a^*, c_1, \ldots, c_m) \ &= \ \tp(a'_0, b'_{\pi(1)}, \ldots, b'_{\pi(m)}, a''_0, b''_{\pi(1)}, \ldots, b''_{\pi(m)}). \nonumber
\end{align}
Then $\crd(a_0) = \{b_1, \ldots, b_m\}$, $\crd(a'_0) = \{b'_1, \ldots, b'_m\}$ and
$\crd(a'') = \crd(a''_0) = \{b''_1, \ldots, b''_m\}$. 
As the equivalence relation $\acl_\mcM(x) = \acl_\mcM(y)$ is assumed to be trivial it follows from
Lemma~\ref{every element is in dcl of its coordinates} that 
$a_0 = a$, $a'_0 = a'$ and $a''_0 = a''$.
Therefore~(\ref{the elements corresponding to a}) implies that
\begin{align*}
&\tp(a, a'') = \tp(a, a^*) \text{ and } \tp(a', a'') = \tp(a', a^*), \text{ and hence} \\
&\tp_\mcM(a, a'') = \tp_\mcM(a, a^*) \text{ and } \tp_\mcM(a', a'') = \tp_\mcM(a', a^*).
\end{align*}
As $\mcM$ is a homogeneous and binary we get
\[ \tp_\mcM(a, a', a'') = \tp_\mcM(a, a', a^*)\]
and hence $\tp(a, a', a'') = \tp(a, a', a^*)$.
Since (by definition) $\crd(a) = \acl(a)  \cap C_1$ and similarly for $a', a''$ and $a^*$, there are 
permutations $\sigma_1, \sigma_2, \sigma_3$ of $\{1, \ldots, m\}$ such that
\begin{align}\label{the long sequences have the same type}
&\tp(b_1, \ldots, b_m, b'_1, \ldots, b'_m, b''_1, \ldots, b''_m) \\ 
= \ 
&\tp(b_{\sigma_1(1)}, \ldots, b_{\sigma_1(m)}, b'_{\sigma_2(1)}, \ldots, b'_{\sigma_2(m)}, c_{\sigma_3(1)}, \ldots, c_{\sigma_3(m)}). 
\nonumber
\end{align}
We have
\[ \{b_1, \ldots, b_m\} \cap \{b'_1, \ldots, b'_m\} \cap \{c_1, \ldots, c_m\} \ \neq \ \es \]
because $b_1 = b'_1 = c_1$ belongs to this intersection, while
\[ \{b_1, \ldots, b_m\} \cap \{b'_1, \ldots, b'_m\} \cap \{b''_1, \ldots, b''_m\} \ = \ \es \]
by the choice of these elements.
This contradicts~(\ref{the long sequences have the same type}).
\hfill $\square$

\end{document}